\documentclass{ws-idaqp}

\usepackage{ifthen}
\newboolean{article}
    \setboolean{article}{true}
\newboolean{report}
\newboolean{book}
\newboolean{letter}
\newboolean{german}
\newboolean{italian}
\newboolean{nobaselinestretch}
   \setboolean{nobaselinestretch}{true}
\newboolean{nosectionappendix}
\newboolean{oldtoc}
\newboolean{nosectionequn}
\newboolean{notheorem}

\ifthenelse{\boolean{german}}{
    \usepackage{german}}{}

\usepackage[latin1]{inputenc}

\usepackage{amsmath}
\usepackage{amssymb}
\usepackage[mathscr]{eucal}

\ifthenelse{\boolean{notheorem}}{}{
    \usepackage{theorem}}

\ifthenelse{\boolean{nobaselinestretch}}{}{
    \renewcommand{\baselinestretch}{1.25}}

\newenvironment{env}[2]{\begin{#1}#2\end{#1}}{}
    \newcommand{\beq}[1]{\begin{env}{equation}{#1}}
    \newcommand{\beqn}[1]{\begin{env}{equation*}{#1}}
    \newcommand{\bal}[1]{\begin{env}{align}{#1}}
    \newcommand{\baln}[1]{\begin{env}{align*}{#1}}
    \newcommand{\bga}[1]{\begin{env}{gather}{#1}}
    \newcommand{\bgan}[1]{\begin{env}{gather*}{#1}}
    \newcommand{\bflal}[1]{\begin{env}{flalign}{#1}}
    \newcommand{\bflaln}[1]{\begin{env}{flalign*}{#1}}
    \newcommand{\bmu}[1]{\begin{env}{multline}{#1}}
    \newcommand{\bmun}[1]{\begin{env}{multline*}{#1}}
    \newcommand{\bsp}[1]{\begin{env}{split}{#1}}

    \newcommand{\eeq}{\end{env}}
    \newcommand{\eeqn}{\end{env}}
    \newcommand{\eal}{\end{env}}
    \newcommand{\ealn}{\end{env}}
    \newcommand{\ega}{\end{env}}
    \newcommand{\egan}{\end{env}}
    \newcommand{\eflal}{\end{env}}
    \newcommand{\eflaln}{\end{env}}
    \newcommand{\emu}{\end{env}}
    \newcommand{\emun}{\end{env}}
    \newcommand{\esp}{\end{env}}

\newcommand{\lf}{\vspace{2ex}}

\renewcommand{\bf}[1]{\textbf{#1}}
\renewcommand{\it}[1]{\textit{#1}}
\renewcommand{\sc}[1]{\textsc{#1}}
\renewcommand{\sf}[1]{\textsf{#1}}
\renewcommand{\sl}[1]{\textsl{#1}}
\renewcommand{\tt}[1]{\texttt{#1}}

\newcommand{\hl}[1]{\bf{\it{#1}}}
\newcommand{\mrm}[1]{\mathrm{#1}}

\newcommand{\msf}[1]{\text{\small$\sf{#1}$}}

\newcommand{\cmc}[1]{\mathcal{#1}}
\newcommand{\eus}[1]{\mathscr{#1}}
\newcommand{\euf}[1]{\mathfrak{#1}}
\newcommand{\bb}[1]{\mathbb{#1}}
\newcommand{\msmall}[1]{{\setlength{\arraycolsep}{.6ex}\text{\small$#1$}}}

\newcommand{\mtiny}[1]{{\setlength{\arraycolsep}{.3ex}\text{\tiny$#1$}}}
\newcommand{\nbd}[1]{$#1$\nobreakdash--}
\newcommand{\ol}[1]{\overline{#1}}

\newcommand{\wh}[1]{\widehat{#1}}

\newcommand{\vt}{\vartheta}

\newcommand{\vp}{\varphi}
\newcommand{\om}{\omega}
\newcommand{\Om}{\Omega}
\newcommand{\abs}[1]{\left\lvert#1\right\rvert}
\newcommand{\norm}[1]{\left\lVert#1\right\rVert}

\newcommand{\bfam}[1]{\bigl(#1\bigr)}
\newcommand{\Bfam}[1]{\Bigl(#1\Bigr)}
\newcommand{\AB}[1]{\langle#1\rangle}

\newcommand{\BAB}[1]{\Bigl\langle#1\Bigr\rangle}
\newcommand{\CB}[1]{\{#1\}}
\newcommand{\bCB}[1]{\bigl\{#1\bigr\}}
\newcommand{\BCB}[1]{\Bigl\{#1\Bigr\}}
\newcommand{\SB}[1]{[#1]}

\newcommand{\Matrix}[1]{\begin{pmatrix}#1\end{pmatrix}}
\newcommand{\SMatrix}[1]{\msmall{\Matrix{#1}}}

\newcommand{\tMatrix}[1]{\mtiny{\Matrix{#1}}}
\newcommand{\rtMatrix}[1]{\raisebox{.3ex}{\tMatrix{#1}}}

\newcommand{\set}[2][]{
    \ifthenelse{\equal{#1}{}}{
        \CB{#2}}{
        \CB{#1~|~#2}}}
\newcommand{\bset}[2][]{
    \ifthenelse{\equal{#1}{}}{
        \bCB{#2}}{
        \bCB{#1~|~#2}}}
\newcommand{\Bset}[2][]{
    \ifthenelse{\equal{#1}{}}{
        \BCB{#2}}{
        \BCB{#1~\big|~#2}}}
\newcommand{\zero}{\CB{0}}

\DeclareMathOperator{\ls}{\normalfont\msf{span}}

\DeclareMathOperator{\id}{\normalfont\msf{id}}

\renewcommand{\dim}{\operatorname{\msf{dim}}}

\newcommand{\C}{\bb{C}}

\newcommand{\E}{\bb{E}}
\newcommand{\V}{\bb{V}}
\newcommand{\bI}{\bb{I}}

\newcommand{\N}{\bb{N}}

\newcommand{\R}{\bb{R}}

\newcommand{\Z}{\bb{Z}}
\newcommand{\cA}{\cmc{A}}
\newcommand{\cB}{\cmc{B}}

\newcommand{\cD}{\cmc{D}}
\newcommand{\cF}{\cmc{F}}

\newcommand{\cI}{\cmc{I}}

\newcommand{\sB}{\eus{B}}

\newcommand{\sE}{\eus{E}}
\newcommand{\sF}{\eus{F}}

\newcommand{\sL}{\eus{L}}

\newcommand{\sN}{\eus{N}}

\newcommand{\G}{\Gamma}

\ifthenelse{\boolean{nosectionequn}}{}{
    \numberwithin{equation}{section}
    \renewcommand{\theequation}{\thesection.\arabic{equation}}}

\ifthenelse{\boolean{article}\or\boolean{letter}\or\boolean{nosectionequn}}{
    \setboolean{nosectionappendix}{true}}{}
\ifthenelse{\boolean{nosectionappendix}}{}{
    \renewcommand{\appendix}{
        \chapter*{\appendixname}
        \addcontentsline{toc}{chapter}{\appendixname}
        \renewcommand{\thesection}{\Alph{section}}
        \setcounter{section}{0}}}
   
\ifthenelse{\boolean{report}\or\boolean{book}}{
    }{}

\ifthenelse{\boolean{notheorem}}{}{
        \newcommand{\notename}{Note.}
        \newcommand{\mnname}{Mathematical note.}
        \newcommand{\enname}{End of the note.}
        \newcommand{\definame}{Definition.}
        \newcommand{\propname}{Proposition.}
        \newcommand{\lemname}{Lemma.}
        \newcommand{\exname}{Example.}
        \newcommand{\exername}{Exercise.}
        \newcommand{\remname}{Remark.}
        \newcommand{\obname}{Observation.}
        \newcommand{\thmname}{Theorem.}
        \newcommand{\corname}{Corollary.}
        \newcommand{\proofname}{Proof.}
    \ifthenelse{\boolean{german}}{
        \renewcommand{\mnname}{Mathematische Notiz.}
        \renewcommand{\enname}{Ende der Notiz.}
        \renewcommand{\exname}{Beispiel.}
        \renewcommand{\exername}{Übung.}
        \renewcommand{\remname}{Bemerkung.}
        \renewcommand{\obname}{Beobachtung.}
        \renewcommand{\thmname}{Satz.}
        \renewcommand{\corname}{Korollar.}
        \renewcommand{\proofname}{Beweis.}}{}
    \ifthenelse{\boolean{italian}}{
        \renewcommand{\mnname}{Nota matematica.}
        \renewcommand{\enname}{Fina della nota.}
        \renewcommand{\definame}{Definizione.}
        \renewcommand{\propname}{Proposizione.}
        \renewcommand{\exname}{Esempio.}
        \renewcommand{\exername}{Esercizio.}
        \renewcommand{\remname}{Nota.}
        \renewcommand{\obname}{Osservazione.}
        \renewcommand{\thmname}{Teorema.}
        \renewcommand{\corname}{Corollario.}
        \renewcommand{\proofname}{Dimostrazione.}

       \renewcommand{\bibname}{Bibliografia}
       
       \renewcommand{\appendixname}{Appendice}
       \renewcommand{\contentsname}{Indice}
       \renewcommand{\listfigurename}{Elenco delle figure}
       \renewcommand{\listtablename}{Elenco delle tabelle}
       
       \renewcommand{\figurename}{Figura}
       \renewcommand{\tablename}{Tabella}

       }{}
    \theoremheaderfont{\normalfont\bfseries}
    \theoremstyle{change}
        \theorembodyfont{\rmfamily}
            \newtheorem{emp}{}[section]
                \newcommand{\bemp}[1][]{
                    \begin{emp}\hskip-\labelsep\bf{#1}\hskip\labelsep}
                \newcommand{\eemp}{\end{emp}}
\newtheorem{itemp}[emp]{}
                \newcommand{\bitemp}[1][]{
                    \begin{itemp}\hskip-\labelsep\bf{#1}\hskip\labelsep\normalfont\itshape}
                \newcommand{\eitemp}{\end{itemp}}
            \newtheorem{note}[emp]{\notename}
                \newcommand{\bnote}{\begin{note}}
                \newcommand{\enote}{\end{note}}
            \newtheorem{mn}[emp]{\mnname}
                \newcommand{\bnm}{\begin{mn}~\begin{quotation}\renewcommand{\baselinestretch}{1}\small\noindent\ignorespaces}
                \newcommand{\enm}{\end{quotation}\hfill\bf{\enname}\end{mn}}
            \newtheorem{ex}[emp]{\exname}
                \newcommand{\bex}{\begin{ex}}
                \newcommand{\eex}{\end{ex}}
            \newtheorem{exer}[emp]{\exername}
                \newcommand{\bexer}{\begin{exer}}
                \newcommand{\eexer}{\end{exer}}
            \newtheorem{defi}[emp]{\definame}
                \newcommand{\bdefi}{\begin{defi}}
                \newcommand{\edefi}{\end{defi}}
            \newtheorem{rem}[emp]{\remname}
                \newcommand{\brem}{\begin{rem}}
                \newcommand{\erem}{\end{rem}}
            \newtheorem{ob}[emp]{\obname}
                \newcommand{\bob}{\begin{ob}}
                \newcommand{\eob}{\end{ob}}
        \theorembodyfont{\normalfont\itshape}
            \newtheorem{thm}[emp]{\thmname}
                \newcommand{\bthm}{\begin{thm}}
                \newcommand{\ethm}{\end{thm}}
            \newtheorem{prop}[emp]{\propname}
                \newcommand{\bprop}{\begin{prop}}
                \newcommand{\eprop}{\end{prop}}
            \newtheorem{cor}[emp]{\corname}
                \newcommand{\bcor}{\begin{cor}}
                \newcommand{\ecor}{\end{cor}}
            \newtheorem{lem}[emp]{\lemname}
                \newcommand{\blem}{\begin{lem}}
                \newcommand{\elem}{\end{lem}}
\newenvironment{empn}[1]{\lf\noindent\bf{#1}\ignorespaces\hskip\labelsep}{\lf}
		\newcommand{\bempn}[1]{\begin{empn}{#1}}
		\newcommand{\eempn}{\end{empn}}
		\newcommand{\bitempn}[1]{\begin{empn}{#1}\normalfont\itshape}
		\newcommand{\eitempn}{\end{empn}}
                \newcommand{\bnmn}{\begin{empn}{\mnname}~\begin{quotation}\renewcommand{\baselinestretch}{1}\small\noindent\ignorespaces}
                \newcommand{\enmn}{\end{quotation}\hfill\bf{\enname}\end{empn}}
		\newcommand{\bexn}{\begin{empn}{\exname}}
		\newcommand{\eexn}{\end{empn}}
		\newcommand{\bexern}{\begin{empn}{\exername}}
		\newcommand{\eexern}{\end{empn}}
		\newcommand{\bdefin}{\begin{empn}{\definame}}
		\newcommand{\edefin}{\end{empn}}
		\newcommand{\bremn}{\begin{empn}{\remname}}
		\newcommand{\eremn}{\end{empn}}
		\newcommand{\bobn}{\begin{empn}{\obname}}
		\newcommand{\eobn}{\end{empn}}
		\newcommand{\bthmn}{\bitempn{\thmname}}
		\newcommand{\ethmn}{\eitempn}

\newcommand{\qedsymbol}{~\hfill$\oblong$}
    \newcounter{proof}[emp]
    \newenvironment{Proof}[1]{
        \vspace{1ex}
        \renewcommand{\item}[1][\stepcounter{proof}(\roman{proof})]%
            {##1\hskip\labelsep}
        \noindent\bf{#1\hskip\labelsep}}{
        \nolinebreak\qedsymbol}
     \renewcommand{\proof}[1][\proofname]{
        \begin{Proof}{#1}\ignorespaces}
    \renewcommand{\qed}{\end{Proof}}
    \newcommand{\noqed}{
        \renewcommand{\qedsymbol}{}
        \end{Proof}}}
    \ifthenelse{\boolean{italian}}{
        \renewcommand{\proofname}{Dimostrazione.}}{}

\newcounter{OP}
		\newcommand{\bOP}{\stepcounter{OP}\begin{empn}{Open Problem \arabic{OP}:}}
		\newcommand{\eOP}{\end{empn}}

\usepackage{lmodern}

\usepackage{stmaryrd}

\usepackage[matrix,arrow,curve]{xy}

\usepackage{hyperref}
\begin{document}

\markboth{Michael Skeide}{Algebraic Central Limit Theorems}

\bibliographystyle{plain}

\title{Algebraic Central Limit Theorems:\\A Personal View on One of Wilhelm's Legacies\\~\\\normalfont To the Memory of Wilhelm von Waldenfels}


\author{Michael Skeide}

\address{Dipartimento di Economia, Universit\`a\ degli Studi del Molise, Via de Sanctis\\
86100 Campobasso, Italy\\
skeide@unimol.it
}

\maketitle

\begin{abstract}
Bringing forward the concept of convergence in moments from classical random variables to quantum  random variables is what leads to what can be called algebraic central limit theorem for (classical and) quantum random variables. I reflect in a very personal way how such an idea is typical for the spirit of doing research in mathematics as I learned it in Wilhelm von Waldenfels's research group in Heidelberg.
\end{abstract}


\ccode{AMS Subject Classification (2020): }



\section{A very personal introduction}

After finishing my physics studies in 1990 with the \it{Diplomarbeit}\cite{Ske90} (master thesis) at the Theoretical Particle Physics in Heidelberg, my supervisor told me that I cannot continue with ``that sort of rigorous mathematics'' at that institute, and sent me off to Wilhelm von Waldenfels -- a move I am eternally grateful for! So, I became Wilhelm's PhD-student.

Apart from many other things (such as an outstandingly good example for what a functioning research group is (see also Footnote \ref{osFN}), or such as total honesty in science%
\footnote{
Honesty, so total that much later I was surprised to have to learn that even in mathematics honesty is not necessarily a granted virtue.}%
), in Wilhelm's research group I learned two things%
\footnote{
These things were to be learned from him in the first place. But I would like to emphasize that also all group members were exemplary.}
I consider determining for the rest of my mathematical life: 
\baln{
\text{\hl{Be elementary!}}
\text{~~~~~~~~~~~~and~~~~~~~~~~~~}
\text{\hl{Be rigorous!}}
}\ealn
Regarding ``Be rigorous!'', after Heidelberg's Theoretical Particle Physics sent me away for being ``too rigorous'', it might come as a surprise that I still had to learn being rigorous better. In fact, I still remember as if it was yesterday, when I was giving my first mathematical talk in Wilhelm's \it{Oberseminar}%
\footnote{ \label{osFN}
``\it{Oberseminar}'' -- the dictionary says ``postgraduate seminar'' (I am not convinced that this translation meets the right spirit) -- is the German terminology for the weekly seminar when the research group -- of course, including students who still have to graduate -- meets to present and discuss the latest news. We would meet \bf{always}, independently on whether there was a talk scheduled or not, and while chatting along with a cup of coffee in our hands, frequently somebody would go to the blackboard and improvise a talk on their latest problems -- and more than once the problems came closer to their solution taking in the feedback from the group.
}
presenting the results of my master thesis, when at a certain point Wilhelm would shout at me ``\it{Schlamperei!}''%
\footnote{
Doing a work in a sloppy, slovenly, messy, untidy, ... way. Not a very friendly comment ...
}%
, after I had messed up. And even after my PhD, when I proposed some nonsense in the seminar, he would say ``You are a doctor of mathematics, you have to know that such a thing cannot work!'' -- adding then a to-the-point explanation, why: ``The concept of Hamel basis does not work together nicely with the topology of a Hilbert space.''

Regarding ``Be elementary!'', this is -- hoping we are being rigorous always -- what this little note is dedicated to.

Being \it{elementary} does not mean being \it{simple}, or \it{easy} (or even \it{trivial}).%
\footnote{
Many referees seem to live under the illusion that an elementary proof it trivial and, thus, not worth of publication.
}
Being as \it{elementary} as possible means to reduce the prerequisites (other, possibly sophisticated theories, theorems and their proofs) of a proof to a minimum. Starting from the premise that a mathematician to fully understand (the proof of) a theorem should understand also the prerequisites (and their proofs), it is clear that being elementary in the proof of a result increases enormously the number of readers who can claim to understand the result. Being \it{elementary} means not to \it{break a butterfly on the wheel}%
\footnote{
In German we say, translated literally, not to \it{shoot with canons at sparrows}.
}%
.

Wilhelm loved elementary proofs (and produced many of them), tailored to deal with a specific problem instead of aiming at maximum generality; but he would not be dogmatic about them. He would not tell his students that a proof of their results \it{has to be} more elementary; but he would recognize an elementary proof if he sees one -- and would get very enthusiastic about it. Just receiving his (rare!) praises on such occasions is enough motivation to wish to get more of them. Let me illustrate how this influenced me with some mathematics:

\bex
Suppose for $-1<q<1$ we have an operator $\alpha\in\sB(H)$ satisfying the \nbd{q^2}CCR
\beqn{
\alpha\alpha^*-q^2\alpha^*\alpha
~=~
1-q^2.
}\eeqn
(Well, it is $a:=\frac{\alpha}{\sqrt{1-q^2}}$ that fulfills the \nbd{q^2}CCR $aa^*-q^2a^*a=1$.) Putting $\G:=1-\alpha^*\alpha$, by easy algebraic computations one verifies $\alpha\G=q^2\G\alpha$. Applying to the \nbd{q^2}CCR the reverse triangle inequality, taking also into account that $\norm{\alpha^*\alpha}=\norm{\alpha\alpha^*}$, we see that $\alpha$, hence, also $\G$ are contractions.

If $q\ne0$, by elementary Hilbert space arguments, the eigenspace of $\G$ to the eigenvalue $0$ is a reducing subspace; we may take it away, that is, we may assume that $\G$ is injective.

To see more clearly how $\alpha$ acts on $H$ (that is, to do the representation theory of the \nbd{q^2}CCR), one may proceed in two different ways:
\begin{itemize}
\item
\begin{itemize}
\item
One may show that the spectrum of $\G$ is $\CB{1,q^2,q^4,\ldots}\cup\zero$.

\item
By applying to the normal operator $\G$ the spectral theorem (and recalling that $\G$ is assumed injective), we get $\G=\sum_{k=0}^\infty E_kq^{2k}$, where $E_k$ is the projection onto the eigenspace of $\G$ to the eigenvalue $q^{2k}$.

\item
By $\alpha\G=q^2\G\alpha$, we see that $\alpha^*$ and $\alpha$ go forth and back between the eigenspaces $E_kH$ and $E_{k+1}H$, and we deduce the form of $\alpha^*$ as weighted right shift.
\end{itemize}
This is more or less how Vaksman and Soibelman\cite{VaSo88} proceeded to find the irreducible representation of the quantum group $SU_q(2)$. Unless you call the spectral theorem (together with all knowledge needed to prove it for bounded operators) elementary, this method is not elementary.

\item
\begin{itemize}
\item
One shows (by a certainly not trivial, but completely elementary proof) the following lemma (see Lemma 2.1\cite{ScSk98}): $\G$ injective $\Longrightarrow$ $\lim_{k\to\infty}\alpha^k=0$ in the strong operator topology.

\item
Since $\G$ is a contraction we may form the Neumann series over $q^{2k}\G$ to find $\frac{1}{1-q^{2k}\G}=\sum_{k=n}^\infty(q^{2k}\G)^n$. We may further define
\beqn{
P_k
~:=~
{\alpha^*}^k\frac{1}{1-q^2\G}\cdots\frac{1}{1-q^{2k}\G}\alpha^k
}\eeqn
($k\in\N$) and $P_0=1$, and put $E_k:=P_k-P_{k+1}$ $(k\in\N_0)$. It is completely elementary to show the $P_k$ is a decreasing projection-valued function that, by the lemma, converges strongly to $0$. Therefore, the $E_k$ determine a spectral measure on $\N_0$. (See Section 2\cite{Ske98b}.)

\item
It is completely elementary to show that $\G=\sum_{k=0}^\infty E_kq^{2k}$ (in norm!) and that $\frac{\alpha^*}{\sqrt{1-q^{2(k+1})}}\upharpoonright E_kH$ and $\frac{\alpha}{\sqrt{1-q^{2(k+1})}}\upharpoonright E_{k+1}H$ define an inverse pair of unitaries between $E_kH$ and $E_{k+1}H$. (See again Section 2\cite{Ske98b}.) From this, the same representation of $\alpha^*$ as weighted right shift follows.
\end{itemize}
This is more or less the representation theory of $SU_q(2)$ from the PhD-thesis\cite{Ske94}. If you agree that the Neumann series and how to deal with commutation relations under the sum is elementary, then so is this proof.
\end{itemize}
(It is noteworthy that the same method in the case $q=0$, where $\alpha^*$ is simply an isometry, looks awkward. On the one hand, the $E_k$ still exist but are much simpler: $E_k={\alpha^*}^k\alpha^k-{\alpha^*}^{k+1}\alpha^{k+1}$. But $\G$ is never injective and the $E_k$ do not correspond to its spectral measure. Still the $E_k$ give back $\alpha^*$ as a (no longer weighted) right shift. The condition that replaces that $\G$ had to be injective, is now one of the many equivalent conditions to say that the isometry $\alpha^*$ is pure or has no unitary part; see, for instance, Section 2\cite{BhSk15}. Recall that $SU_0(2)$ is no longer a quantum group, but only a quantum monoid.)
\eex

\lf
Wilhelm wrote towards the end of the abstract of Ref.\cite{Wal73}: ``The mathematical techniques are mainly the \it{well-known methods of nearly elementary probability theory}.'' It is a pity that the expression ``nearly elementary probability'' did not enter standard terminology for presenting probability.%
\footnote{
This brings to my mind the book\cite{Let93} entitled ``Probabilità Elementare'' which I used, following an advice from Franco Fagnola, to prepare myself for the competition yielding me my present job. I am really wondering if the students of Computer Science with whom, as the back-flap says, this textbook has been tested, would agree to call it very elementary. A book that really sticks to finitely additive probabilities for quite a time before introducing measure theory, is Parthasarathy's textbook\cite{Par05}. In my lecture notes\cite{Ske14s} (also for Computer Science) I stay with finite additivity to the very end. (Of course, the probabilities I consider \bf{are} \nbd{\sigma}additive; but to work at this level, one need not know that. The discussion is independent on which integration theory you put as a basis.)
}
One (out of many) convincing example(s) of being \it{nearly elementary} in Wilhelm's work is surely everything that has to do with algebraic central limit theorems; this is what I try to illustrate in Section \ref{CLTsec} (the Giri-von Waldenfels algebraic CLT) and \ref{bCLTsec} (the Speicher-von Waldenfels general algebraic CLT), while in Section \ref{FCLTsec} I discuss some of its impact, presenting in the appendix the CLT for boolean conditional independence.

\section{CLT \it{versus} algebraic CLT} \label{CLTsec}

The classical C(entral) L(imit) T(heorem), in its simplest form for an i(ndependent) i(dentically) d(distributed) sequence of real random variables $(X_n)$ on a probability space,  asserts: If $X_n\sim X$ with $\E X=0$ and $\V X=1$, then the sequence of random variables
\beqn{
S_N
~:=~
\frac{X_1+\ldots+X_N}{\sqrt{N}}
}\eeqn
converges in distribution to a random variable $Y$ (not necessarily on the same probability space) with standard normal distribution $N(0,1)$, that is,
\beqn{
\lim_{N\to\infty}\E f(S_N)
~=~
\E f(Y)
}\eeqn
for every bounded measurable function $f$ on $\R$.

There are several subclasses of the bounded measurable functions such that it is sufficient to check convergence for the subclass only:
\begin{itemize}
\item
$f=\bI_I$ being indicator functions of intervals $I$. (That is, $\E\bI_I(X)=P(X\in I)=:\mu(I)$ with $\mu$ being the probability measure of $X$, explaining the name convergence in distribution.)

\item
The continuous functions with bounded support. (More or less the Riesz representation theorem.)

\item
Passing to complex valued functions, we have the collection of functions $f_\om:=e^{i\om\bullet}$ ($\om\in\R$). The function $\om\mapsto\wh{\mu}(\om):=\E f_\om(X)=\E e^{i\om X}$ is the \it{characteristic function} or \it{Fourier transform} of $\mu$, $\wh{\mu}$.
\end{itemize}
The simplest proof of the above CLT is surely by characteristic functions.%
\footnote{
It takes advantage of the fact that, if $X\sim\mu$ and $Y\sim \nu$ are independent, then $X+Y\sim\mu*\nu$ and $\wh{\mu*\nu}=\wh{\mu}\,\wh{\nu}$.
}
That indicator functions suffice to have convergence in distribution, is a statement that makes part of Lebesgue integration; it is up to you if you consider it elementary. But none of the proofs that checking convergence on the other two classes suffices, is exactly elementary.%
\footnote{
The proof of the Riesz theorem as presented in Section 50 of the book\cite{RiSz-Na82} by Riesz and Sz.-Nagy, is a valid predecessor of the proof of the spectral theorem for bounded self-adjoint operators presented in Section 107 of the same book. The proof that pointwise convergence of the characteristic functions is enough, uses concepts such as \it{tightness} of a sequence of probability measures Banach-Alaoglu-type arguments.
}
So accepting that proving pointwise convergence of the characteristic functions is elementary, the proof that this suffices for having the CLT is not.

But apart from these technicalities, the problem when we wish to pass from classical probability to quantum (or noncommutative, or algebraic) probability, and try to have a CLT also there, is a different one: To say what convergence in distribution is for random variables, one has to say what a function of a random variable is (and compute their expectations). (Real)  quantum variables are (self-adjoint) elements in a unital (no longer necessarily commutative) \nbd{*}algebra $\cA$ and the expectation functional is replaced by a functional $\vp\colon\cA\rightarrow\C$ which is linear,  positive (that is, $\vp(a^*a)\ge0$), and normalized (that is, $\vp(1)=1$), that is, $\vp$ is a \it{state}. To say what for a function $f$ on $\R$ and for a self-adjoint element $a\in\cA$ we understand by $f(a)$, we need what is called a \it{functional calculus} at $a$: For each $a=a^*\in\cA$ a unital homomorphism $F_a$ from a unital algebra of functions on $\R$ containing $\id_\R$ into $\cA$ such that $F_a(\id_\R)=a$.

The only functional calculus which we have for elements in a general (algebraic) quantum probability spaces $(\cA,\vp)$ is for polynomials -- and every more powerful calculus (in necessarily more restricted classes of quantum probability spaces) starts from the polynomial calculus and extends it (as far as the subclass allows). The  best known instances are holomorphic functional calculus for Banach \nbd{*}algebras, continuous functional calculus for \nbd{C^*}algebras, and measurable functional calculus for von Neumann algebras. The latter two are forms of the spectral theorem and also work for not necessarily bounded self-adjoint (in the functional analytic definition of adjoint) operators on a Hilbert space $H$.%
\footnote{ \label{vNFN}
The latter is used to consider included classical probability into von Neumann probability spaces $(\cA,\vp)$ where $\cA\subset\sB(H)$ is a von Neumann algebra and $\vp$ a normal state on $\cA$. Looking, then, at (not necessarily bounded) operators on $H$ that are \it{affiliated} with $\cA$ (roughly speaking, operators that have a polar decomposition where the partially isometric factor is in $\cA$ and where the  positive radial part has a spectral measure in $\cA$ allowing to define all sorts of bounded functions of it in $\cA$), one gets enough unbounded operators but everything is determined by what $\vp$ does on projections. Gleason type theorems tell when a function defined on the lattice of projections is the restriction of a normal state.

Given an algebraic probability space $(\cA,\vp)$, we get by GNS-construction a pre-Hilbert space $\cD=\cA/\sN$, where $\sN$ is the space of length-zero-elements of the form $(a,b):=\vp(a^*b)$, with a representation $\pi$ of $\cA$ by adjointable operators on $\cD$ and cyclic vector $\Om:=1+\sN$ (that is, $\pi(\cA)\Om=\cD$) such that $\AB{\Om,\pi(a)\Om}=\vp(a)$. Note that the algebraic domain $\cD$ is invariant under all $\pi(a)$. If a self-adjoint element $a$ has a $\pi(a)$ that is essentially self-adjoint on $\cD$, then we get its spectral measure. If this is true for all self-adjoint elements of $\cA$, then we may pass to the von Neumann algebra on $H=\ol{\cD}$ generated by the projections of all these spectral measures, and for all self-adjoint elements $a\in\cA$ the operator $\pi(a)$ is affiliated with this von Neumann algebra.

If the GNS-representation $\pi$ is faithful, we may identify $\cA$ as algebra of operators on $\cD$. Sufficient for that is, if $\vp$ is \it{faithful}, that is, if $\vp(a^*a)=0$ implies $a=0$; but this is by far not necessary. (Any vector state on $\sB(H)$ ($\dim H\ge2$) is a counter example.) It is sufficient and necessary, if $\vp$ is a trace -- as in classical probability where the algebra is commutative.
}
But, if any functional calculus starts with polynomials, then why not looking at convergence at polynomials to begin with? Why not simply look at convergence of the moments?
 
Well, one main answer to this question is that, in classical probability, not all random variables $X$ possess all moments $\E X^m$ ($m\in \N_0$). But for the classical CLT stated in the beginning of this section, existence of the second moment is enough. (There are even better versions without the latter condition.) The trick is that for convergence in distribution we need to know only the expectations of the random variables $f(X)$ for bounded functions,  which all exist \it{a priori} and computing them does not require to know the expectations of the random variables $X^m$. (This is similar to the procedure in Footnote \ref{vNFN}, reducing everything to those projections that are indicator function of the (possibly) unbounded quantum random variables in question.)

In an algebraic quantum probability space, it is part of the definition that the state $\vp$ exists on the whole algebra -- and this includes, for $a=a^*\in\cA$ the elements $a^m$. But, if all moments exist, then checking convergence of moments is enough also classically provideded the moments determine the limit distribution: If $Y,X_1,X_2,\ldots$ are random variables all possessing all (finite!) moments and if there is no other distribution than that of $Y$ having the same moments $\E Y^m$, then convergence of $\E X_n^m\xrightarrow{~n\to\infty~}\E Y^m$ for any $m$ implies $X_n\xrightarrow{~n\to\infty~}Y$ in distribution.%
\footnote{
See, for instance, Theorem 11.4.1 in the nice book\cite{Ros00} by Rosenthal. Unfortunately, this is not a standard topic in a first course in measure theoretic probability; many -- otherwise, really excellent -- books do not discuss it.
}
This is what we take as notion of convergence for algebraic central limit theorem:

A sequence $a_n$ of self-adjoint elements in a quantum probability space $(\cA,\vp)$ \hl{converges in moments} to a self-adjoint element $b$ in a quantum probability space $(\cB,\psi)$ if
\beqn{
\lim_{n\to\infty}\vp(a_n^m)
~=~
\psi(b^m)
}\eeqn
for all $m\in\N_0$. (We entirely ignore the question whether these are moments of a classical probability distribution, or whether they determine it.)

With this, after specifying what is \it{independent}, we are able to formulate the analogue of the classical CTL from above. Suppose $a=a^*$ is in the quantum probability space $(\cA_0,\vp_0)$. Define $\cA$ to be the subalgebra of $\bigotimes_{n\in\N}\cA_0\ni a_1\otimes a_2\otimes\ldots$ generated by elements
\beqn{
a_1\otimes\ldots\otimes a_N\otimes 1\otimes\ldots~~~~~~(N\in\N;a_i\in\cA_0),
}\eeqn
and (well-)define the linear functional $\vp$ on $\cA$ by
\beqn{
\vp(a_1\otimes\ldots\otimes a_N\otimes 1\otimes\ldots)
~:=~
\vp_0(a_1)\ldots\vp_0(a_N);
}\eeqn
obviously, $\vp$ is a state on $\cA$.%
\footnote{
$\cA$ is the inductive limit of the unital \nbd{*}algebras $\cA^N:=\bigotimes_{n=1}^N\cA_0$ under the canonical embeddings of the $\cA^N$ into the first $N$ factors of $\cA^{N+M}$, and $\vp$ is the linear functional defined on this inductive limit by the compatible family $\vp^N:=\vp_0^{\otimes N}$ of functionals on $\cA^N$.
}
Let us denote the \nbd{n}th factor $1\otimes\ldots\otimes 1\otimes\cA_0\otimes 1\otimes\ldots$ in $\cA$ by $\cA_n$ and let us identify $a_n\in\cA_0$ with the corresponding element in $\cA_n$. Especially, for $a\in\cA_0$ denote by $a_n$ the copy of $a$ in $\cA_n$. Define
\beqn{
s_N
~:=~
\frac{a_1+\ldots+a_N}{\sqrt{N}}.
}\eeqn

\bthm \label{aCLT}
If $a=a^*\in\cA_0$ with $\vp_0(a)=0$  and $\vp_0(a^2)=1$, then
\beqn{
\lim_{N\to\infty}\vp(s_N^m)
~=~
\frac{1}{\sqrt{2\pi}}\int_{-\infty}^\infty x^me^{-\frac{x^2}{2}}\,dx,
}\eeqn
the \nbd{m}th moment of the normal law $N(0,1)$, for all $m\in\N_0$. In other words, the sequence $s_N$ converges in moments to a standard gau{\ss}ian (quantum) random variable.
\ethm

This is a (quite) special case of the algebraic CLT by Giri and von Waldenfels (Theorem 1 in Ref.\ \refcite{GiWa78}). We discuss on the proof, under much more general circumstances, in the next section.

\bnote
 Giri and von Waldenfels\cite{GiWa78} is probably the first time that an algebraic quantum CLT has been proved (that is, by convergence of moments requiring not more than existence of all moments of $a$). There is the earlier result by Cushen and Hudson\cite{CuHu71} -- a contribution of another great man, Robin Hudson, to be commemorated in this volume -- where the proof is by considerably more analytic (pushing forward the characteristic function approach, leading to distributions on the Weyl algebra), but limited to sequences of pairs $(p_n,q_n)$ of quantum variables that fulfill \it{a priori} the Heisenberg commutation relations $p_nq_n-q_np_n=\frac{1}{i}$. Of course, the Heisenberg commutation relations limits what $p_n$ and $q_n$ can be. (For instance, neither of them can be bounded in the GNS-representation; also, under regularity conditions on the GNS-representation -- that is, on the state -- the representing operators are quite determined by the Stone- von Neumann theorem.) It is noteworthy that the general version of Theorem \ref{aCLT} (Theorem 1\cite{GiWa78}) plus Theorem 2\cite{GiWa78} allows to recover Theorem 1\cite{CuHu71} in this algebraic way. (See Example 3\cite{GiWa78} to Theorem 1\cite{GiWa78}.) Also a version for Fermion commutation relations can be shown. (See von Waldenfels\cite{Wal78}.)
\enote

\lf
Theorem \ref{aCLT} and the classical version we spoke about in this section, are \it{individual} central limit theorems in the sense that there is a single distribution (of $X$ in the classical case, and of $a$ in $(\cA_0,\vp_0)$ in the quantum case) and the convergence is to a distribution of a single random variable (of $Y$ in the classical case, and of $b$ in $(\cB,\psi)$ in the quantum case). Already in the Hudson-Cushen CLT mentioned in the preceding note, we are speaking about a sequence of pairs with a fixed joint distribution. The Giri-von Waldenfels CLT, of which we reported in Theorem \ref{aCLT} only the one-variable-version, is actually a multivariate CLT: It is a CLT for sequences of \bf{families} of (quantum) random variables that  are identically distributed  (that is, they have identical joint distributions) -- and that are independent in a certain sense, namely, \it{tensor independent} (explaining the occurrence of the infinite tensor product occurring in the construction of the quantum probability space $(\cA,\vp)$ in Theorem \ref{aCLT}).

The multivariate case and situations different from tensor independence, is what leads us right into the next section. We conclude the present section with some more remarks.

\bnote
Independence -- or \it{quantum} independence -- is one of the topics which Wilhelm's work helped crucially unfolding. \it{Tensor independence}, the most classical one (as it captures classical independence when applied to commutative quantum probability spaces), is considered by many the ``most natural'' one, in the sense that when physicists wish to model two quantum physical systems as subsystems of a single one in such a way that one has nothing to do with the other, they take tensor products. But -- unlike in classical  probability -- there are (sometimes, many) more independences in quantum probability. One of the most famous is Voiculescu's\cite{Voi85} \it{free independence} (or \it{freeness}); it is the ``most noncommutative'' one and in a very strict (functorial) view point is may be considered the ``only noncommutative'' one (it is the only \it{unital} independence that survives the transition from states to conditional expectations; see Theorem 2.5\cite{Ske04} and see also Section \ref{FCLTsec} of the present notes). A (non-unital) independence that is even more noncommutative, is \it{monotone independence} (Lu\cite{Lu97}, Muraki\cite{MurN97}); here we have that ``$A$ independent of $B$'' does not imply ``$B$ independent of $A$''.  \it{Boolean independence} (Speicher and Woroudi\cite{SpWo97}) can be traced back to von Waldenfels\cite{Wal73}. There is \it{\nbd{q}independences} (Bozejko and Speicher\cite{BoSp91a}) for $q\in\SB{-1,1}$, interpolating tensor (or Boson) independence ($q=1$) and Fermion independence\cite{Wal78} ($q=-1$) going over freeness ($q=0$).

Starting (probably) with Accardi, Schürmann, and von Waldenfels\cite{ASW88}, one of Wilhelm's students, Michael Schürmann, has taken on the challenge to examine independence axiomatically\cite{MSchue95b} and to examine the quantum Lévy processes associated with each such independence. For tensor independence, the research culminated conclusively in his habilitation thesis\cite{MSchue93}. This work also discusses \nbd{q}independence --  but this does not belong to the axiomatic independences. In fact, axiomatic independences are very few; see Speicher\cite{Spe97}, Ben Ghorbal and Schürmann\cite{BGSc02}, and Muraki\cite{MurN02}. On the other hand, dropping the one or the other requirement of axiomatic independence, there can be very many -- in fact, enough so, to derive, in the single-random-variable-case any central distribution as central limit distribution for a suitable notion of independence. (See Theorem 6.1(v) in Accardi and Bozejko\cite{AcBo98}.)

We do not discuss this in detail; there will be more about independences and Lévy processes in Michael Schürmann's contribution to this volume. We mentioned independence only to underline how much it means to say that the CLT in the next section's title covers the technical part of CLTs for \bf{every} independence.
\enote


\section{The Speicher-von Waldenfels general algebraic CLT} \label{bCLTsec}

Let us right start with the theorem from Speicher and von Waldenfels\cite{SpWa94} -- and then explain it.

\bitemp[Central Limit Theorem\cite{SpWa94}.] \label{SWCLT}
Let $(\cA,\vp)$ be a quantum probability space. Let $J$ be a fixed index set. Consider elements $b_n^{(j)}\in\cA$ $(j\in J,n\in\N)$ such that for each $j\in J$ there is $j'$ with $(b_n^{(j)})^*=b_n^{(j')}$, which fulfill the following assumptions.

\begin{itemize}
\item[$i)$]
We have
\beqn{
\vp\Bfam{b_{\sigma(1)}^{(j_1)}\ldots b_{\sigma(n)}^{(j_n)}}
~=~0
}\eeqn
for all $n\in\N$, all $(j_1,\ldots,j_n)\in J^n$ and all $\sigma\colon\{1,\ldots,n\}\rightarrow\N$ with the property that there exists $k\in\N$ such that $\#\sigma^{-1}(k)=1$.

\item[$ii)$]
For all $n\in\N$ there exists a constant $C_n<\infty$, such that
\beqn{
\abs{\vp\Bfam{b_{\sigma(1)}^{(j_1)}\ldots b_{\sigma(n)}^{(j_n)}}}
~\le~
C_n
}\eeqn
for all $(j_1,\ldots,j_n)$ $\in J^n$ and all $(\sigma(1),\ldots,\sigma(n))\in\N^n$.

\item[$iii)$]
We have an invariance of all second order correlations under order preserving injections, i.e.
\beq{ \label{oiv}
\vp\Bfam{b_{\vartheta(\sigma(1))}^{(j_1)}\ldots b_{\vartheta(\sigma(n))}^{(j_n)}}
~=~
\vp\Bfam{b_{\sigma(1)}^{(j_1)}\ldots b_{\sigma(n)}^{(j_n)}}
}\eeq
for all $n\in\N$, all $(j_1,\ldots,j_n)\in J^n$, all $\sigma\colon\{1,\ldots,n\}\rightarrow\rm{Im}(\sigma)$ with the property that for all $k\in\rm{Im}(\sigma)$ we have $\#\sigma^{-1}(k)=2$, and for all order preserving injective mappings $\vartheta\colon\rm{Im}(\sigma)\rightarrow\N$.
\end{itemize}

For each $N\in\N$ we define
\beqn{
S_N^{(j)}
~:=~
\frac{b_1^{(j)}+\cdots+b_N^{(j)}}{\sqrt{N}}.
}\eeqn
Then we have for all $n\in\N$ and all $(j_1,\ldots,j_n)\in J^n$
\beqn{
\lim_{N\to\infty}\vp\bfam{S_N^{(j_1)}\ldots S_N^{(j_n)}}
~=~
0,
}\eeqn
if $n$ is odd, and
\beqn{
\lim_{N\to\infty}\vp\bfam{S_N^{(j_1)}\ldots S_N^{(j_n)}}
~=~
\frac{1}{\bigl(\frac{n}{2}\bigr)!}\sum_{\substack{\sigma\colon\CB{1,\ldots,n}\rightarrow\CB{1,\ldots,\frac{n}{2}}\\\forall k\in\CB{1,\ldots,\frac{n}{2}}\colon\#\sigma^{-1}(k)=2}}\vp\Bfam{b_{\sigma(1)}^{(j_1)}\ldots b_{\sigma(n)}^{(j_n)}},
}\eeqn
if $n$ is even.
\eitemp

How to interpret this?
\begin{itemize}
\item
The index $j\in J$ labels a selection of $\#J$ elements $b^{(j)}$ in a unital \nbd{*}algebra $\cA_0$.

In fact, if we take the free unital algebra generated by $\#J$ symbols $b^{(j)}$ and define an involution by ${b^{(j)}}^*:=b^{(j'(j))}$ (the $j'$ that exists by the hypotheses of the theorem for any $j$; it need not be unique, but any fixed choice works). Then for each $n\in\N$ we can define a homomorphism $\alpha_n\colon\cA_0\rightarrow\cA$ sending $b^{(j)}$ to $b^{(j)}_n$.

We continue denoting $a_n:=\alpha_n(a)$ for $a\in\cA_0$, as in Theorem \ref{aCLT}. $S^{(j)}_N$ is what you usually would write down in a multivariate CLT.

\item
The condition in (ii) is the only technical condition that in the end yields convergence. For concrete independences, this is, usually, fulfilled automatically.

The condition in (i) means for $n=1$ that the $b^{(j)}$ are sent to central elements in $\alpha_k(\cA_0)$. ($\leadsto$ \bf{Central} Limit Theorem.) For bigger $n$ it yields what is called frequently the \it{singleton condition}: If in a monomial of generators from several $\alpha_i(\cA_0)$, there is exactly one coming from a certain $\alpha_k(\cA_0)$, then $\vp$ at that monomial is $0$. (This is a condition any independence fulfills. It means that in a sum over monomials of degree $n$, only those contribute where, for each $k$, there occurs either no factor or there occur at least two. This is responsible for the fact that in a CLT the normalization $\frac{1}{N}$ (yielding the \it{law of large numbers}) may be replaced with the weaker normalization $\frac{1}{\sqrt{N}}$ without producing divergence.)

The result means not only that the moments on the left-hand side converge. It means that in the sum that occurs when expanding the product of $n$ factors of the type $b_1^{(j)}+\cdots+b_N^{(j)}$ only those monomials contribute to the limit where from the generators of each $\alpha_k(\cA_0)$ there are either no factors or exactly two factors. Under these general conditions this is the best one can do to compute the limit distribution. How to compute each number under the sum is what the specific independence in question tells us; see the Examples \ref{cCLTex} below and the appendix. 

\item
If, in (iii), we require \eqref{oiv} for all $\sigma\colon\CB{1,\ldots,n}\rightarrow\N$, we get what is called a sequence of quantum probability subspaces $\alpha_n(\cA_0)$ with \it{spreadable} (joint) distribution. (For most independences, if the $\alpha_n(\cA_0)$ are independent and identically distributed (see below), then their joint distribution is spreadable.) Fixing, under spreadability, to $\sigma$ that is constant $n$ (only the $n$th copy of $\cA_0$, $\alpha_n(\cA_0)$, is involved) we get that
\beqn{
\vp_0(a)
~:=~
\vp(a_n)
}\eeqn
does not depend on $n$, hence, well-defines a state $\vp_0$ on $\cA_0$ turning it into a quantum probability space $(\cA_0,\vp_0)$: The $\alpha_n(\cA_0)$ are identically distributed.%
\footnote{
This is, why in a consequent dualization of ``probability space'' into ``quantum probability space'', quantum random variables are homomorphisms from a \nbd{*}algebra $\cA_0$ into a quantum probability space. (This dualizes classical random variables, which are functions \bf{from} a probability space \bf{into} a measurable space.) The fixed \it{quantum measurable space} is $\cA_0$, and the the quantum random variables are the $\alpha_n$. So, that $\vp\circ\alpha_n$ does not depends on $n$ is what it means to say the $\alpha_n$ are \it{identically distributed}.
}
The condition in (iii) alone does not guarantee identically distributed but only \it{identically correlated}: The correlations $\vp(b_n^{(j)}b_n^{(\ell)})$ for all $j,\ell\in J$ do not depend on $n$.
\end{itemize}

\noindent
In the limit formula for moments of even degree $n$ the sum is taken over functions $\sigma\colon\CB{1,\ldots,n}\rightarrow\CB{1,\ldots,\frac{n}{2}}$ such that for each $k$ in the codomain $\CB{1,\ldots,\frac{n}{2}}$ there are exactly two $i_k\ne j_k$ in the domain $\CB{1,\ldots,n}$ such that $\sigma(i_k)=\sigma(j_k)=k$. So, such a function determines a partition of the domain $\CB{1,\ldots,n}$ into pairs $\CB{i_k,j_k}$, that is, a \it{pair partition}. Two such functions $\sigma,\sigma'$ determine the same pair partition if and only if there is a permutation $\pi\in S_{\frac{n}{2}}$ of the codomain$\CB{1,\ldots,\frac{n}{2}}$ such that $\sigma'=\pi\circ\sigma$. The map $\pi$ considered as map into $\N$ is injective, yes, but it is order preserving if and only if $\pi=\id$. So, under the condition in (iii) we have no invariance under $\pi$. However, \bf{if} we have invariance under all $\pi$, equivalently, if the invariance under $\vt$ is not only for order preserving injections but for all injections, then all $\pi\circ\sigma$ give the same contribution. Since there are $(\frac{n}{2})!$ of them, we may cancel out the factor in front of the sum, and take the sum over all pair partitions. The occurrence of pair partitions in central limit theorems and that frequently for determining the moments of the central limit distribution one has ``just'' to count the number of certain (sub)classes of pair partitions, plays an important role in the research of another of Wilhelm's students, Roland Speicher. 

\bex \label{cCLTex}
We discuss some one-variable-cases, that is, $\#J=1$ and we write $b^{(j)}=b$. (Note that in each of the following examples the growth condition in (ii) is fulfilled with $C_n=1$.)
\begin{enumerate}
\item \label{cex1}
In the case of Theorem \ref{aCLT}, tensor independence, we have $\vp_0(b^2)=1$, and for the functions $\sigma$ occurring in the sum we have
\beqn{
\vp(b_{\sigma(1)}\ldots b_{\sigma(n)})
~=~
\vp(b_1^2)\ldots\vp(b_{\frac{n}{2}}^2)
~=~
\vp_0(b^2)^{\frac{n}{2}}
~=~
1.
}\eeqn
There are $\binom{n}{2}$ possibilities to select the pair that is sent to $1$, then, $\binom{n-2}{2}$ possibilities to select the pair that is sent to $2$, and so forth. We get that there are
\beqn{
\binom{n}{2}\binom{n-2}{2}\ldots\binom{2}{2}
~=~
\frac{n!}{2^{\frac{n}{2}}}
}\eeqn
functions $\sigma$. Consequently,
\beqn{
\lim_{N\to\infty}\phi(S_N^n)
~=~
\frac{n!}{2^{\frac{n}{2}}(\frac{n}{2})!}
~=~
n!!
}\eeqn
(for even $n$). These are both the number of pair partitions and the even moments of the normal law $N(0,1)$.

\item \label{cex2}
For free independence, we have that
\beqn{
\vp(a^{(1)}_{\sigma(1)}\ldots a^{(n)}_{\sigma(n)})
~=~
0
}\eeqn
for any selection of central $a^{(k)}\in\cA_0$ ($1\le k\le n$) and any $\sigma$ such that $\sigma(k-1)\ne\sigma(k)\ne\sigma(k+1)$ ($1<k<n$).%
\footnote{
Well, writing $\sigma(k)\ne\sigma(k+1)$ ($1\le k<n$) would be enough, and we will do so in the sequel. But here we wish to emphasize on the position $k$ and that is is surrounded by elements from different subalgebras, and for that we do accept a certain redundancy.
}
Writing a general $a\in\cA_0$ as $a=\vp(a)-(a-\vp(a))$ (where $a-\vp(a)$ is central), this becomes a recursion to compute $\vp$ on any monomial in terms of $\vp_0$ alone.

If, in the sum of the CLT, the pair partition of $\sigma$ has no pair consisting of neighbouring elements (that is, it $\sigma(k-1)\ne\sigma(k)\ne\sigma(k+1)$), then for the central elements $b$ we have $\vp(b_{\sigma(1)}\ldots b_{\sigma(n)})=0$. If there is a pair $\CB{i_k,j_k}$ with (without loss of generality) $j_k=i_k+1$, then, rewrite
\beqn{
b_{i_k}b_{j_k}
~=~
\vp(b_{i_k}b_{j_k})+(b_{i_k}b_{j_k}-\vp(b_{i_k}b_{j_k})).
}\eeqn
The term arising from the summand $\vp(b_{i_k}b_{j_k})$ has the same structure but with degree reduced from $n$ to $n-2$ (and a factor $\vp(b_{i_k}b_{j_k})=1$ in front of it). The central element $b_{i_k}b_{j_k}-\vp(b_{i_k}b_{j_k})$ remaining in the term arising from the second summand, is a singleton -- and none of the consecutive operations to reduce the degree further by means of the recursion can change this; its contribution is, therefore, $0$. We see:
\begin{itemize}
\item
 Only pair partitions that can be reduced by consecutive elimination of next neighbour pairs contribute to the sum.
 
\item
The contribution of each such so-called \it{noncrossing pair partition} (for the way they can be represented graphically) is $1$.
\end{itemize}
To determine the limit of the moments for the free CLT, we are, thus, reduced to the problem to count the number of the noncrossing pair partitions of $\CB{1,\ldots,n}$ ($n$ even). They are known to be given by the \it{Catalan numbers}
\beqn{
\mrm{Catalan}_{\frac{n}{2}}
~:=~
\frac{n!}{\bfam{\frac{n+2}{2}}!\bfam{\frac{n}{2}}!},
}\eeqn
and they coincide with the even moments of the \it{Wigner semicircle law} with density $x\mapsto\frac{\sqrt{2^2-x^2}}{2\pi}$ on $\SB{-2,2}$.

\item \label{cex3}
For boolean independence, we have that
\beqn{
\vp(a^{(1)}_{\sigma(1)}\ldots a^{(n)}_{\sigma(n)})
~=~
\vp(a^{(1)}_{\sigma(1)})\ldots\vp(a^{(n)}_{\sigma(n)})
}\eeqn
for \bf{any} (that is, not necessarily central) selection of $a^{(k)}\in\cA_0$ ($1\le k\le n$) and any $\sigma$ such that $\sigma(k-1)\ne\sigma(k)\ne\sigma(k+1)$ ($1<k<n$). We see that in the sum over $\vp(b_{\sigma(1)}\ldots b_{\sigma(n)})$ only those $\sigma$ contribute that have a pair partition with only next-neighbour pairs (each of them with a factor $\vp_0(b^2)^{\frac{n}{2}}=1$). There is only one such pair partition. The even moments of the boolean CLT limit distribution are just $1$, the distribution, therefore, given by the probability measure $\frac{\delta_{-1}+\delta_1}{2}$.
\end{enumerate}
\eex

\noindent
Examples where the invariance condition in (iii) really only holds for order preserving injections $\vt$, are the \nbd{q}CLT (Schürmann\cite{MSchue91c}) and the CLT for monotone independence (with the \it{arcsine law} as limit distribution).

\section{Fock, conditional expectation, and module,  and all that} \label{FCLTsec}

A generalization that was surely beyond what Wilhelm might have dreamed of is the passage from \it{independence} in a state to \it{conditional independence} in a \it{conditional expectation} (also known as \it{independence with amalgamation} or as \it{operator-valued independence}). In classical probability, where everything is done with emphasis on the probability measure(s) and the expectation functional $\E$ is a second-order object, \it{conditional expectation} -- a map on algebras of random variables -- is, actually, much closer to its generalization to quantum probability: A positive idempotent \bf{onto} a subalgebra that is also a bimodule map for the bimodule structure of the big algebra over its subalgebra.%
\footnote{
In classical probability, one requires, additionally, that the measure (that is, the expectation functional) is preserved. This makes, given the subalgebra, a conditional expectation unique, while existence is a consequence of the Radon-Nikodym theorem (only the version for probability measures). Several authors do require the same in quantum probability: Their conditional expectation on a (usually, von Neumann) quantum probability space should preserve the state. (Also this forces uniqueness, while existence is characterized by Takesaki's famous theorem that the state preserving conditional expectation exists if and only if the \it{modular automorphism group} of the (faithful normal) state leaves the subalgebra (globally) invariant.)
}
Note that we do not even require that the big (unital) \nbd{*}algebra $\cA$ is a quantum probability space. (No state!) We just speak about a \it{\nbd{\cB}quantum probability space} $(\cA,\Phi)$ where $\Phi$ is a conditional expectation onto a unital \nbd{*}subalgebra $\cB$ of $\cA$. (Typically, it is required  that the unit of $\cB$ coincides with the unit of $\cA$.%
\footnote{
There are conditional expectations onto subalgebras that have a different unit. Very typical is $a\mapsto pap$ where the subalgebra $p\cA p$ of $\cA$ has its own unit $p$. This happens, for instance, in monotone and boolean independence, where (at least some of) the independent subalgebras have different units but are, usually, \it{expected} (that is, the image of a conditional expectation). But for our new noncommutative scalars $\cB$, we expect that their unit acts as unit on the whole thing.
}%
) As already observed in Skeide\cite{Ske99a}, Theorem \ref{SWCLT} remains true in this setting without changing a word in its proof.

Among many others, the axiomatic independences have been promoted to conditional independences in a \nbd{\cB}quantum probability space $(\cA,\Phi)$. There is \it{free conditional independence} (or \it{freeness with amalgamation}) also in Voiculescu\cite{Voi85}, which underlies Roland Speicher's habilitation thesis\cite{Spe98}. There is monotone conditional independence in Skeide\cite{Ske04} and Popa\cite{MPop08}.%
\footnote{
We take the occasion to mention that in Skeide\cite{Ske04} there is a flaw (regarding the non-valid attempt to identify the $\cB$ of the nonunital subalgebra with that of the containing algebra. But the moment defining formula is there -- and Popa's\cite{MPop08} coincides with that.
}
Boolean conditional independence has been formulated in Skeide\cite{Ske00}. In the free case, there is an individual CLT by Voiculescu\cite{Voi95} and a multivariate version by Speicher\cite{Spe98}. The monotone CLT is due to Popa\cite{MPop08}. A multivariate CLT for the boolean case is presented in the appendix of the present notes.

To cite p.71 of Lance\cite{Lan95}: ``The enquiring reader will already have asked the question:'' What happened to tensor independence? Well, the answer (Skeide\cite{Ske04}) is that it disappeared -- at least in full generality there is no such thing as the tensor product of algebras \bf{over} a common subalgebra $\cB$, and writing down the moment defining formula simply leads to ill-defined nonsense. (In a certain sense, this is nice: There is only one ``truly noncommutative'' (namely, over noncommutative scalars $\cB$!) axiomatic conditional independence that is unital: Freeness! See also Roland Speicher's contribution to this volume.) A way out is to limit conditional tensor independence to so-called \it{centered} \nbd{\cB}quantum probability spaces, that is, the algebra is generated by elements that commute with the subalgebra $\cB$, and limit the moment defining formula to these (from where they extend well-definedly to the whole). This has been introduced (as \it{Bose \nbd{\cB}independence}) in Skeide\cite{Ske03a} (preprint 1996) and Skeide\cite{Ske99a} discusses the CLT.%
\footnote{
A warning: Being centered is a tough requirement for a bimodule. (It means it is a subbimodule of a free bimodule, generated by a vector space. Only \nbd{\sB(H)}bimodules fulfill this automatically; see Skeide\cite{Ske98}.) And, on the other hand, in this framework all desirable commutation relations can be imposed in the centered case, by imposing them on elements of $\cB$. (See Skeide\cite{Ske99a} for \nbd{q}commutation relations.) As always, conditional independence is the more interesting the ``smaller'' $\cB$ is (also in the sense the ``lesser'' structure $\cB$ has) as compared with $\cA\supset\cB$. More interesting is when the relevant \nbd{\cB}bimodules are generated by elements that fulfill (nontrivial) commutation relations with elements of $\cB$; see, for instance, the treatment of the \it{square of white noise} in Accardi and Skeide\cite{AcSk00a}.

Conditional tensor independence does have advantages in classical probability (where it always exists).
}

The proof of CLTs is done most conveniently by -- simultaneously -- identifying the limit moments as moments of certain operators on \it{Fock-type} objects. (This is true already for scalar independence, but the bigger effect is obtained in the conditional versions.) The fact that we replace the scalars $\C$ by noncommutative (that is, not necessarily commutative!) scalars $\cB$ cries for Hilbert modules -- in the same way \nbd{*}algebras with positive functionals cry for Hilbert spaces. For states on a unital \nbd{*}algebra we have the well-known GNS-construction; see Footnote \ref{vNFN}. For CP-maps -- and a conditional expectation is a CP-map! -- we have Paschke's GNS- construction\cite{Pas73}. We do not discuss this for general CP-maps, but directly in the case of conditional expectations, which is a good deal more similar to the GNS-construction for states. (See the discussion in Section 4.4\cite{Ske01}.) If $\Phi\colon\cA\rightarrow\cB\subset\cA$ is a conditional expectation, then define on $\cA$ the semiinner product $\AB{a,a'}:=\Phi(a^*a')$, divide out the length-zero elements $\sN$ and complete the quotient $\cA/\sN$ (a pre-Hilbert \nbd{\cB}module) to obtain a Hilbert \nbd{\cB}module $\sE$. Then $a(a'+\sN):=aa'+\sN$ (well-)defines a left action of $\cA$ on $\sE$ (we say, $\sE$ is turned into a \it{correspondence} from $\cA$ to $\cB$), and with the cyclic vector $\om:=1+\sN$ (in fact, $\ol{\cA\om}=\sE$) we get
\beqn{
\AB{\om,a\om}
~=~
\Phi(a).
}\eeqn
Note that $\om$ is in the center of $\sE$ ($b\om=\om b$), so that $\sE\supset\om\cB=\cB$.

In fact, Speicher\cite{Spe98} identified the moments of the CLT-distribution in conditional freeness%
\footnote{
Note that this means freeness (or free independence) in a conditional expectation. It has nothing to do with the ``conditionally free'' random variables in Bozejko \it{et al.}\cite{BLS96}, which is ``scalar business''.
}
as moments of certain sums $a^*(x)+a(y)$ of the creation and annihilation operators on the \it{full Fock module}.%
\footnote{ \label{FockFN}
The first mentioning of the full Fock module might have been in Accardi and Lu\cite{AcLu96}, but (looking also at the preprint data) is was almost simultaneous with Pimsner\cite{Pim97}and Speicher\cite{Spe98}. The time was simply ripe for that.

It is noteworthy that Accardi and Lu\cite{AcLu96}, motivated by the concrete physical example they study in that paper, pass immediately from full Fock modules to so-called \it{interacting Fock modules}. The so-called \it{weak coupling limit} in a conditional expectation  they compute may actually be considered as a sort of CLT in its own right. Only in Skeide\cite{Ske98} it has been observed that (putting the right left action) the module carrying the limit distribution is, actually, a full Fock module. Note, however, that this does not mean that the limiting objects would be conditionally free in some sense; the occurring creators and annihilators associated to certain disjoint time intervals do not have orthogonal arguments. (See also Note \ref{Findn}.)

As a spin-off of the paper\cite{AcLu96}, there is the notion of \it{interacting Fock space} (IFS), so only the scalar version, introduced and studied systematically by Accardi \it{et al.}\cite{ALV97,Lu97,AcBo98,AcSk08,GeSk20b}: A pre-Hilbert space (allowing technically to deal with unbounded operators without problems) $\cI=\bigoplus_{n\in\N_0}H_n$ for a family of pre-Hilbert spaces $H_n$ with $H_0=\C\Om$ and a linear map $a^*\colon H\rightarrow\sL(\cI)$ ($H$ another pre-Hilbert space) such that
\beqn{
\ls a^*(H)H_n
~=~
H_{n+1}.
}\eeqn
(See Refs.\ \refcite{AcSk08,GeSk20b} for this way to define IFS.) If we assume that $\cI$ is \it{adjointable} in the sense that every $a^*(f)$ has an adjoint $a(f)\colon\cI\rightarrow\cI$, then one would view $a^*(x)+a(x)$ and $a^*(y)+a(y)$ as independent in the vacuum state $\vp:=\AB{\Om,\bullet\Om}$ (inducing, thus, a notion of independence) whenever $\AB{x,y}=0$. We mentioned already that using IFS, every central distribution occurs as CLT-distribution.
}
Popa\cite{MPop08} did the same for monotone conditional independence. The appendix of the present notes discusses the boolean case. Before passing to the appendix, we conclude the discussion with the following considerations.

\bnote \label{Findn}
The observation that all classical (\nbd{\R}valued) Lévy processes may be represented in a very canonical way on the symmetric Fock space $\G(L^2(\R_+,K))$ as a sum of suitable creation, conservation, and annihilation operators, goes probably back to Parthasarathy and Schmidt\cite{PaSchm72}. This includes, of course, that every single random variable with infinitely divisible distribution may be represented on the Fock space. But more interesting is the interplay with independence in this representation. Roughly, if we look how the increment $X_t-X_s$ from $s<t$ to $t$ of the Lévy process are represented, then it turns out that the representing operators ``live on $L^2(\SB{s,t},K)\subset L^2(\R_+,K)$'', meaning that they are in
\beqn{
\sL\bfam{\G(L^2(\SB{s,t},K))}\otimes\id
~\subset~
\sL\Bfam{\G(L^2(\SB{s,t},K))\otimes\G(L^2(\SB{s,t},K)^\perp)}
}\eeqn
in the well-known factorization $\G(H_1\oplus  H_2)$``$=$''$\G(H_1)\otimes\G(H_2)$ of the symmetric Fock space. (We write $\sL$ because the operators are unbounded.) This means, the increments, looked at as quantum random variables, are tensor independent in the \it{vaccum state} $\vp=\AB{\Om,\bullet\Om}$, where $\Om$ denotes the \it{Fock vacuum}. (Schürmann has generalized this to all quantum Lévy processes.)

``All'' Hilbert spaces (infinite-dimensional and separable) are isomorphic -- and representing a single distribution by opertors on a Fock-type space is almost as uninformative: We know, it is always possible. This becomes more interesting, when the question is whether the quantum random variables represented by creation and annihilation (and, sometimes, conservation) operators to disjoint intervalls are indedent in some sense.

Of course, the more structure the Fock-type space has ($=$the less general it is!), the more informative is also what we get for a single variable. (For instance, the Fock space has a ``strong'' inherent structure, if all distributions of $a(f)^*+a(f)$ ($a^*(f)$ some creation operator and its adjoint $a(f)$) are from the same restricted convolution family (for some sort of independence) of distributions. In fact, in the symmetric Fock space we get only normal laws $N(0,\sigma^2)$ and in the full Fock space we get only (scaled) semicircular laws. To include, on the symmetric Fock space, general Lèvy processes (such as Poisson processes) creators and annihilators already are not enough: One also needs conservation operators.

We see (see also Footnote \ref{FockFN}) that representing a certain distribution on a Fock space (or module) is one thing. The question whether several of them done on the same Fock space (or module) might be (conditionally) independent is a totally different one.
\enote

\section*{Appendix: Conditional boolean CLT}

Let $\cB$ be a unital \nbd{C^*}algebra, and let $E$ be a \nbd{\cB}correspondence. The \it{boolean Fock module} over $E$ is the \nbd{\cB}correspondence $\cF_B(E):=\cB\oplus E$.%
\footnote{
This is an interacting Fock \it{module} $\cI=\bigoplus_{n\in\N_0}E_n$ (with $E_1=E$ and $E_n=\zero$ for $n\ge2$) in its own right.
}
We denote by $\om:=1\in\cB$ the \it{vacuum}. For $x\in E$ the creator $a^*(x)$ is defined as
\beqn{
a^*(x)
~:=~
\SMatrix{0&0\\x&0}
~\in~
\sB^a(\sF_B(E))
~=~
\SMatrix{\sB^a(\cB)&\sB^a(E,\cB)\\\sB^a(\cB,E)&\sB^a(E)}
~=~
\SMatrix{\cB&E^*\\E&\sB^a(E)}
}\eeqn
(for the last equality recall that $\cB$ is unital), with adjoint $a(x):=a^*(x)^*=\rtMatrix{0&x^*\\0&0}$.

Recall that for a \nbd{\cB}quantum probability space $(\cA_0,\Phi_0)$ we have the GNS-correspondence $\sE$ with the cyclic and \bf{central} (the latter, because the CP-map $\Phi_0$ is a conditional expectation) unit vector $\om$ such that $\Phi_0=\AB{\om,\bullet\om}$. This decomposes $\sE$ into
\beqn{
\sE
~=~
\cB\oplus E
}\eeqn
where the summand isomorphic to $\cB$ is $\om\om^*\sE=\om\cB$ and where $E=\om^\perp$. For any $a\in\cA_0$ we denote by
\beqn{
\SMatrix{\alpha&\beta^*\\\gamma&\delta}
}\eeqn
($\alpha$ in $\cB$, $\beta$ and $\gamma$ in $E$, $\delta\in\sB^a(E)$) its action on $\cB\oplus E$.

Recall from Skeide\cite{Ske00} that subalgebras $\cA_k$ of a \nbd{\cB}quantum probability space $(\cA,\Phi)$  are \it{conditionally boolean independent} if
\beqn{
\Phi(a_1\ldots a_n)
~=~
\Phi(a_1)\ldots\Phi(a_n)
}\eeqn
for any $a_k\in\cA_{\sigma(k)}$ ($1\le k\le n$) whenever $\sigma$ is such that $\sigma(k)\ne\sigma(k+1)$ ($1\le k<n$).

\bthmn
Let $(\cA_0,\Phi_0)$ and $(\cA,\Phi)$ be \nbd{\cB}quantum probability spaces, and for each $n\in\N$ let $\pi_n\colon\cA_0\rightarrow\cA$ be a(n almost never unital, so not \nbd{\cB}preserving) homomorphism respecting the bimodule action of $\cB$ (that is, $\pi_n(bab')=ba_nb'$, where, as usual, $a_n:=\pi_n(a)$) such that the subalgebras $\cA_n$ are conditionally boolean independent in $\cA$.

For every selection of elements $b^{(j)}\in\cA_0$ ($j$ in an index set $J$) with $\Phi_0(b^{(j)})=0$, putting
\beqn{
S_N^{(j)}
~:=~
\frac{b^{(j)}_1+\ldots+b^{(j)}_N}{\sqrt{N}},
}\eeqn
we obtain
\beqn{
\lim_{N\to\infty}\Phi(S_N^{(j_1)}\ldots S_N^{(j_n)})
~=~
\BAB{\om,\bfam{a^*(\gamma^{(j_1)})+a(\beta^{(j_1)})}\ldots\bfam{a^*(\gamma^{(j_n)})+a(\beta^{(j_n)})}\om}
}\eeqn
where for each $j$ we denote by $\rtMatrix{\alpha^{(j)}&{\beta^{(j)}}^*\\\gamma^{(j)}&\delta^{(j)}}$ the action of $b^{(j)}$ on $\cB\oplus E$.
\ethmn

\proof
For each limit is sufficient to limit the discussion to the finite subset $\AB{j_1,\ldots,j_n}$ of those $j$ that actually occur. On this, by $\norm{\Phi}=1$, the growth condition in (ii) is, obviously, fulfilled with $C_n=\bfam{\max_k\norm{b^{(j_k)}}}^n$.

As for the scalar case in Example \ref{cCLTex}\eqref{cex3}, we argue that only $\sigma$ with a partition into next neighbour pairs contributes and that the order of the number is $\CB{1,\ldots,\frac{n}{2}}$ that label these pairs do not matter. We find that the only contribution from the sum in the formula in Theorem \ref{SWCLT} that remains is
\beqn{
\Phi_0(b^{(j_1)}b^{(j_2)})\ldots\Phi_0(b^{(j_{n-1})}b^{(j_n)}).
}\eeqn
Taking into account that $\Phi_0(b^{(j)})=0$, hence, $\alpha^{(j)}=0$, we get $\Phi_0(b^{(j)}b^{(j')})=\AB{\beta^{(j)},\gamma^{(j')}}$. So,
\beqn{
\lim_{N\to\infty}\Phi(S_N^{(j_1)}\ldots S_N^{(j_n)})
~=~
\AB{\beta^{(j_1)},\gamma^{(j_2)}}\ldots\AB{\beta^{(j_{n-1})},\gamma^{(j_n)}}
}\eeqn
That this coincides with $\BAB{\om,\bfam{a^*(\gamma^{(j_1)})+a(\beta^{(j_1)})}\ldots\bfam{a^*(\gamma^{(j_n)})+a(\beta^{(j_n)})}\om}$, follows easily from $a(x)\om=0$, from $a^*(x)a^*(y)=0$, and from $a(x)a^*(y)\om=\om\AB{x,y}$.\qed

\newpage

\vspace{1ex}\noindent
\bf{Acknowledgments.~}
Apart from thanking for useful comments by Malte Gerhold and Orr Shalit, and bibliographic assistance from my predecessors in Heidelberg, Michael Schürmann and Roland Speicher, it is my heartfelt wish to thank Wilhelm for his teaching of the heart of how to look at mathematics as I received it through him and his students. Lat, but certainly not least, I wish to thank the referee for an extremely meticulous reading of this manuscript.

\vspace{-2ex}
\newcommand{\Swap}[2]{#2#1}\newcommand{\Sort}[1]{}


\end{document}